\documentclass[12pt]{article}

\newenvironment{Item}[1]{\par #1 }{\par}

\newtheorem{thm}{Theorem}[section]

\newtheorem{prop}[thm]{Proposition}

\newtheorem{corl}[thm]{Corollary}

\newtheorem{rem}[thm]{Remark}

\begin{document}

\title{Maximal ideal space of a commutative coefficient algebra}
\author{ B.K. Kwa\'sniewski, A.V. Lebedev\\
\\
Institute of Mathematics, Bia{\l}ystok University,\\ ul. Akademicka 2,
PL-15-424, Bia{\l}ystok, Poland} \maketitle

\begin{abstract}
The basic notion of the article  is a pair $(\mathcal{A},U)$, where
$\mathcal{A}$ is a commutative $C^*$-algebra and $U$ is a partial isometry
such that $\mathcal{A} \ni a \to UaU^*$ is an endomorphism of
$\mathcal{A}$ and  $U^*U\in \mathcal{A}'$.
We give a description of the maximal ideal space
 of the smallest coefficient $C^*$-algebra
$\overline{E_*(\mathcal{A})}$ of the algebra $C^*(\mathcal{A},U)$
generated by the system $(\mathcal{A},U)$.
\end{abstract}

{\small {\ KEY WORDS: commutative $C^*$-algebra, endomorphism, partial isometry, coefficient algebra }

\section{Introduction. Extensions of $C^*$-algebras by partial isometries
and coefficient algebras}

The notion of a coefficient algebra was introduced  in \cite{lebiedodzij}
in connection with the study of extensions of $C^*$-algebras by partial isometries.

Namely in \cite{lebiedodzij} the authors  investigated the following object.

 Let $H$ be a Hilbert space and  ${\cal A}\subset L(H)$  be a certain $^*$-algebra
 containing the identity $1$ of $L(H)$, the paper was devoted to the description
 of the $C^*-$extensions of $\cal A$ associated with the mappings
\begin{equation}
\label{delta}
\delta(x)=UxU^*,\qquad \delta_*(x)=U^*xU,\qquad x\in L(H)
\end{equation}
where $U\in L(H),\ U\neq 0$. It is clear that
$\delta$ and $\delta_*$ are linear and continuous
($\Vert\delta\Vert =\Vert{\delta_*}\Vert =\Vert U^2\Vert$) maps of $L(H)$ and
$\delta(x^*)=\delta(x)^*$, \ \ $\delta_*(x^*)=\delta_*(x)^*$. When using the powers $\delta^k$ and $\delta_*^k,\ \
k=0,1,2,\ldots $ we
 assume for
convenience that $\delta^0 (x) = \delta_*^0 (x) = x $.

Observe that if $\delta : {\cal A} \to L(H)$ is a morphism then we have
 $$
UU^* = \delta (1) = \delta (1^2) = \delta^2 (1) = (UU^*)^2
  $$
and therefore $U$ is a partial isometry.

In \cite{lebiedodzij} the authors  studied the $C^*$-algebra
$C^*({\cal A},U)$ generated by $\cal A$ and $U$ assuming additionally
that $\cal A$ is the \textbf{coefficient algebra} of $C^*({\cal A},U)$,
by this they  meant  that  $\cal A$ possessed the following three
properties
\begin{equation}
\label{b1}
 {\cal A}\ni a\rightarrow\delta(a)=UaU^*\in {\cal A},
\end{equation}
\begin{equation}
\label{b2}
 {\cal A}\ni a\rightarrow\delta_*(a)=U^*aU\in {\cal A},
\end{equation}
\begin{equation}
\label{b3}
 Ua=\delta(a)U, \qquad a\in {\cal A}
\end{equation}
 Algebras possessing these properties really play the role of
'coefficients' in $C^*({\cal A},U)$ which was shown in \cite{lebiedodzij}, Proposition 2.3.
  telling us  that

{\em if  a $^*$-algebra ${\cal A}$ and $U$ satisfy  conditions (\ref{b1}), (\ref{b2}),
(\ref{b3}) then  the vector space consisting of  finite sums
\begin{equation}
x=U^*a_{\overline{n}}+...+U^*a_{\overline{1}}+a_0+a_1U+...+a_N
U^N,
\end{equation}
where $a_k$, $a_{\overline{l}}\in {\cal A}$ and
$N\in{\bf N}\cup \{0\}$, is a uniformly dense $^*$-subalgebra of
the $C^*$-algebra $C^*({\cal A},U)$.}

It is worth mentioning that  property (\ref{b3}) can
also be written in a different equivalent form
which is stated in the next
\begin{prop}
\label{b4}
{\rm (\cite{lebiedodzij}, Proposition 2.2.)}\ \
Let $\mathcal{A}$ is a $C^*$-algebra of $L(H), \ 1\in \mathcal{A}$ and  $U\in L(H)$
then the following conditions are equivalent
\begin{itemize}
\item[(i)] $Ua=\delta(a)U,\qquad a\in\mathcal{A}$,
\item[(ii)] $U$ is a partial isometry  and
\begin{equation}
\label{A'}
U^*U\in \mathcal{A}'
\end{equation}
where
$\mathcal{A}'$ is the commutant of $\mathcal{A}$,
\item[(iii)]  $U^*U\in \mathcal{A}'$ and $\delta : \mathcal{A} \to \delta (\mathcal{A})$
is a morphism.
\end{itemize}
\end{prop}

Thus property  (\ref{b3}) implies that the mapping
$\delta(a)=UaU^*$,  $a \in{\cal A}$, is a morphism. Hence,
every coefficient algebra is associated with an endomorphism
$\delta$ generated by a partial isometry $U$.

 In \cite{lebiedodzij} the authors gave
 a construction of an  algebra satisfying (\ref{b1}), (\ref{b2}) and (\ref{b3})
starting from an initial algebra that satisfies only some of these conditions
or even does not satisfy any of them. Hereafter we present a part of this construction
for the case of a commutative algebra $\cal A$.

 Let us denote by
$$\overline{E_*({\cal A})}=\overline{\{
\bigcup_{n=0}^\infty\delta_*^n({\cal A})\}}
$$ the  $C^*$-algebra
generated by $\bigcup_{n=0}^\infty\delta_*^n({\cal A})$. It was proved in
\cite{lebiedodzij} (Proposition 4.1.)  that the following statement is true.

\begin{prop}
Let ${\cal A}$ be a commutative $C^*$-subalgebra of  $L(H)$
containing $1$. Let  $\delta$  be an endomorphism of
${\cal A}$ and let $U^*U\in {\cal A}'$ then the
$C^*$-algebra $\overline{E_*({\cal A})}=\overline{\{
\bigcup_{n=0}^\infty\delta_*^n({\cal A})\}}$  is the minimal commutative
  coefficient algebra for $C^*({\cal A},U)$ and both the mappings
$\delta:\overline{E_*({\cal A})}\rightarrow
\overline{E_*({\cal A})}$ and
$\delta_*:\overline{E_*({\cal A})}\rightarrow
\overline{E_*({\cal A})}$ are endomorphisms.
\end{prop}

This proposition leads to the natural problem of obtaining the
description of the maximal ideal space of the coefficient algebra
$\overline{E_*({\cal A})}$ in terms of the maximal ideal space of
$\cal A$ and the action of $\delta$.
Precisely the solution to this problem is the theme of the present  article.
\\

We have to mention that some {\em concrete}  examples of the
description of the maximal ideal space of $\overline{E_*({\cal
A})}$ in the situation when ${\cal A} = C[a,b]$ and $\delta$ is
generated by a continuous mapping $\alpha : [a,b] \to [a,b]$ of a
{\em special} form are given in
\cite{Popovich-Maistrenko}.
\\

The paper is organized as follows. In the second section
we introduce a number of  necessary for our future goals notions and notation
and present some
(mostly known) facts on the structure of endomorphisms of commutative algebras.
Our main result --- the description of the maximal ideal space of
$\overline{E_*({\cal A})}$ is given in section 3.

\section{Endomorphisms of commutative $C^*$-algebras}

 Our starting objects are   a
commutative $C^*$-algebra $\cal A$ containing an
identity $1$ and an endomorphism
\begin{equation}
\label{b0'}
\delta : {\cal A} \to {\cal A} .
\end{equation}

Our first observation is that any
endomorphism $\delta$ generates  a continuous \emph{partial mapping}
on the maximal ideal space
$ M=M(\mathcal{A})$ of  $\mathcal{A}$
which is stated in theorem
\ref{tu1.0}.

{\bf Remark.} This theorem
is a particular case of the general result on the description of
endomorphisms of semisimple Banach
algebras \cite{Lo} and we present its proof for the sake of completness.\\

Since by the Gelfand-Naimark theorem the Gelfand transform defines
an isomorphism $\mathcal{A}\cong C(M)$  we shall identify $\cal A$
and $C(M)$ in all  further considerations.
\begin{thm}
\label{tu1.0} Let  $\cal A$ be a commutative  $C^*$-algebra. Let
$1\in \cal A$ and  $\delta$ be an endomorphism of  $\cal A$. Then
there exists a subset  $\Delta$ of the maximal ideal space  $M$
such that
\begin{description}
\item[i)] the set  $\Delta$ is open and closed,
\item[ii)] the endomorphism  $\delta$ is given by the formula
\begin{equation}
\label{e3.0}
(\delta f)(x)=\left\{ \begin{array}{ll} f(\alpha(x))& ,\ \ x\in
\Delta\\
0 & ,\ \ x\notin \Delta \end{array}\right. ,
\end{equation}
where  $f\in C(M)$ and  $\alpha:\Delta \rightarrow M$ is a continuous mapping.
\end{description}
\end{thm}
\noindent {\bf Proof.}
Let $\Delta\subseteq M$
be the set given by the condition
\begin{equation}\label{Delta}
\tau\in\Delta\Longleftrightarrow \tau(\delta(1))=1.
\end{equation}
First let us observe the closedness and openness of  $\Delta$.
Note that  $\delta(1)=\delta^2(1)$ and therefore the function
$\delta(1)\in C(M)$ is an idempotent,
so its values are either $0$ or $1$.
This implies the closedness and openness of $\Delta$.\\

In terms of the Gelfand transform $a \to \widehat{a}$ we can define the
action of  $\delta$ on  $C(M)$ by means of the formula
\begin{equation} (\delta
\widehat{a})(\tau)=\tau(\delta(a))=\widehat{a}(\delta^*(\tau)),\qquad \tau
\in M,
\end{equation}
where  $\delta^*:A^*\rightarrow A^*$ is
the adjoint  operator to $\delta$. Clearly
\begin{equation}\label{kij}
\delta^*(\tau)=\tau\circ\delta \end{equation}
 is a multiplicative functional and
 $\delta^*(\tau)(1)=\tau(\delta(1))$.
 By   the definition of the set  $\Delta$ we have
\begin{equation}
\tau\notin \Delta \Rightarrow \delta^*({\tau})\equiv 0,
\end{equation}
 \begin{equation}
 \tau\in\Delta \Rightarrow \delta^*(\tau)\in M.
 \end{equation}
Now  defining the map
$\alpha:\Delta\longrightarrow M$
as the restriction of the map  $\delta^*$
\begin{equation}\label{alfaokr}
\alpha=\delta^{*} |_{_\Delta}
\end{equation}
we obtain the statement of the theorem.
\begin{thm}
\label{tu12}
Let the conditions of theorem
\ref{tu1.0} be satisfied and let the mapping  $\alpha:\Delta\rightarrow M$
be given by (\ref{e3.0}). Then
\begin{description}
\item[i)] if \ \ $\mathrm{ker}\delta=\{0\}$\ \ then  $\alpha:\Delta \rightarrow M$
is a surjection,
\item[ii)] if \ \ $\delta(1)=1$\ \ then  $\Delta=M$.
\end{description}
\end{thm}
\noindent {\bf Proof.}
Let
 $\mathrm{ker}\delta=\{0\}$. Then
$\delta$ is an injection and there exists a left inverse to it
$\varrho:\delta({\cal A}) \rightarrow \cal A$. Therefore we have
\begin{equation}
\varrho(\delta(a))=a,\qquad a\in {\cal A}.
\end{equation}
For any  $\tau \in M$ the functional  $\tau\circ \varrho$ being defined on $\delta({\cal A})$ is
nonzero and multiplicative. Therefore there exists its extension
 $\tau_1\in M$ on  $\cal A$ (see
\cite{Dix}, 2.10.2.) and thus
\begin{equation}
\tau_1\circ\delta=\tau.
\end{equation}
Now (\ref{kij}) and  (\ref{alfaokr}) imply the surjectivity of  $\alpha$.\\

Now let  $\delta(1)=1$. This equality means that for any  $\tau \in M$
we have
$\tau(\delta(1))=1$. So  (\ref{Delta}) implies  $\Delta=M$.
The proof is complete.\\

Now we return to the initial object of the article. To this end  $\cal A$
is a  $C^*$-subalgebra of the algebra  $L(H)$ containing the identity operator
$1$ and  $U\in L(H)$ is an operator  such that the mapping
\begin{equation}
\label{b0}
\delta(a)=UaU^*, \ \ a\in \mathcal{A}
\end{equation}
is an endomorphisms of $\cal A$
(this implies in particular that $U$ is a partial isometry).

Note that applying endomorphism $\delta$
$n$ times one  obtains
$$
U^nU^{*n}= \delta^n (1) =  \delta^n (1^2) = (\delta^n (1))^2 = (U^nU^{*n})^2 .
$$
Which means that  $U^n$ is a partial isometry and therefore
$U$ is a power isometry. \\

In view of the form  (\ref{b0}) of the endomorphism  $\delta$ and
relation (\ref{Delta}) it is possible to rewrite theorem
\ref{tu1.0} for the objects considered in the following form
\begin{thm}
\label{tu1.1}
Let  ${\cal A}\subset L(H)$ be a commutative  $C^*$-algebra containing the identity operator
$1$ and let the mapping  $\delta(a)=UaU^*$ be an endomorphism of
the algebra  $\cal A$. Then
\begin{description}
\item[i)] The set  $\Delta=\{\tau\in M: \tau(UU^*)=1\}$ is open and closed,
\item[ii)] On the maximal ideal space $M$ the endomorphism
$\delta$  is given by the formula
\begin{equation}
\label{b-4}
(\delta f)(x)=\left\{ \begin{array}{ll} f(\alpha(x))& ,\ \ x\in
\Delta\\
0 & ,\ \ x\notin \Delta \end{array} \right.
\end{equation}
where  $f\in C(M)$, and  $\alpha:\Delta \rightarrow M$ is a continuous mapping.
\end{description}
\end{thm}

 Moreover theorem  \ref{tu12} implies
\begin{thm}
\label{takietamtwierdzenie}
Let the conditions of theorem  \ref{tu1.1}
be satisfied and  $\alpha:\Delta\rightarrow M$ be given by (\ref{b-4}).
Then we have
\begin{description}
\item[i)] if \,  $U$ is a unitary operator then  $\Delta=M$ and
$\alpha:M \rightarrow M$ is surjective,
\item[ii)] if \,  $U$ is an isometry then
$\alpha:\Delta \rightarrow M$ is surjective,
\item[iii)] if \,  $U^*$ is an isometry then  $\Delta=M$.
\end{description}
\end{thm}

As the next theorem shows the situation  described in theorem
\ref{tu1.1} simplifies considerably if  along with endomorphism
 $\delta$ (\ref{b0}) the mapping
\begin{equation}
\label{delta_*}
\delta_* (a)=U^*aU \ ,\qquad a\in \cal A.
\end{equation}
is an endomorphism of $\cal A$ as well.
\begin{thm}
\label{tu1}
Let  $\cal A$ be a commutative  $C^*$-subalgebra of  $L(H)$
containing the identity
 $1$ and the mappings  $\delta$, $\delta_*$  given by formulae
 (\ref{b0}), (\ref{delta_*}) be endomorphisms of the algebra  $\cal A$.
 Let  $M$ be the maximal ideal space of  $\cal A$. Then
\begin{description}
\item[i)] the sets  $\Delta_1=\{\tau\in M: \tau(UU^*)=1\}$\ \
and  \ \ $\Delta_{-1}=\{\tau\in M: \tau(U^*U)=1\}$  are open and closed,
\item[ii)] in terms of the algebra  $C(M)$ the endomorphism $\delta$
is given by the formula
\begin{equation}
(\delta f)(x)=\left\{ \begin{array}{ll} f(\alpha(x))& ,\ \ x\in
\Delta_1\\
0 & ,\ \ x\notin \Delta_1 \end{array} \right.
\end{equation}
where  $f\in C(M)$ and  $\alpha:\Delta_1 \rightarrow \Delta_{-1}$
is a homeomorphism.
\item[iii)] the endomorphism  $\delta_*$ is given by the formula
\begin{equation}
(\delta_* f)(x)=\left \{ \begin{array}{ll} f(\alpha^{-1}(x))&
,\ \ x\in
\Delta_{-1}\\
0 & ,\ \ x\notin \Delta_{-1} \end{array} \right.
\end{equation}
where  $f\in C(M)$.
\end{description}
\end{thm}
\noindent {\bf Proof.} By theorem  \ref{tu1.1} there exist
closed and open sets  $\Delta_1,\, \Delta_{-1} \subseteq M$
given by the relations
\begin{equation}
\tau\in\Delta_1\Longleftrightarrow \tau(UU^*)=1,
\end{equation}
\begin{equation}
\tau\in\Delta_{-1}\Longleftrightarrow \tau(U^*U)=1
\end{equation}
and continuous mappings  $\alpha:\Delta_1\rightarrow M$,
$\alpha':\Delta_{-1}\rightarrow M$ for which  $\delta$ and
 $\delta_*$ satisfy  (\ref{b-4}) respectively.
 To finish the proof it is enough to show that
$\alpha'=\alpha^{-1}$. This follows from the  relations
\begin{equation}
\tau\in\Delta,\, a\in A\Longrightarrow
\tau(\delta(\delta_*(a))=\tau(UU^*)\tau(a)\tau(UU^*)=\tau(a)
\end{equation}
\begin{equation}
\tau\in\Delta_{-1},\, a\in A\Longrightarrow
\tau(\delta_*(\delta(a))=\tau(U^*U)\tau(a)\tau(U^*U)=\tau(a)
\end{equation}
Or  $(\alpha'\circ \alpha)(\tau)=(\alpha\circ
\alpha')(\tau)=\tau.$ The proof is complete.
\ \\

In the next section we shall make use  of the sets introduced hereafter.

Let
 \begin{equation}
\label{D_n}
\Delta_n=\alpha^{-n}(M), \ \ n=0,1,2, ...
\end{equation}
be  the set on which $\alpha^n$
is defined and let
\begin{equation}
\label{D-n}
\Delta_{-n}=\alpha^{n}(\Delta_n), \ \ n=1,2, ...
\end{equation}
be  the image of
$\alpha^n$.

We have that
\begin{equation}
\label{b-2}
 \alpha^n:\Delta_n\rightarrow  \Delta_{-n},
\end{equation}
\begin{equation}\label{b-3}
\alpha^n (\alpha^m(x)) = \alpha^{n+m}(x),\qquad x\in \Delta_{n+m}.
\end{equation}

\noindent In terms of the multiplicative functionals the sets
$\Delta_n$ can be defined in the following form: for $n>0$
\begin{equation}
\label{e1.0}
\tau\in \Delta_{n} \Longleftrightarrow \forall_{0<k\leq n}\,
\tau(U^kU^{*k})=1,
\end{equation}
\begin{equation}\label{e2.0}
\tau\in \Delta_{-n} \Longleftrightarrow
\exists_{\tau_n\in\Delta_n}\,
 \tau_n\circ\delta^n=\tau .
\end{equation}
Note however that the sequence of projections  $U^kU^{*k}$ is decreasing
and therefore  $\tau(U^nU^{*n})=1$ implies
$\tau(U^kU^{*k})=1$ for  $k<n$. So we can rewrite  condition   (\ref{e1.0})
in the form
\begin{equation}
\label{e1.01}
\tau\in \Delta_{n} \Longleftrightarrow
 \tau(U^nU^{*n})=1.
\end{equation}

\begin{rem}\emph{ In the situation considered in  theorem \ref{tu1}
one can describe the sets  $\Delta_{-n}$ as the sets on which the mapping
$\alpha^{-n}$ is defined.
\\ Moreover in terms of the maximal ideals we have
\begin{equation}
\tau \in \Delta_n \Longleftrightarrow \tau(U^nU^{*n})=1,
\end{equation}
\begin{equation}
\tau \in \Delta_{-n} \Longleftrightarrow \tau(U^{*n}U^n)=1,
\end{equation}
where $n\geq 0$
}
\end{rem}

             \section{Maximal ideal space of a commutative coefficient algebra}

Throughout this section we fix a commutative $C^*$-subalgebra
${\cal A}\subset L(H)$, $1\in {\cal A}$, and a
partial isometry $U\in L(H)$  such that the mapping (\ref{b0}) is an
endomorphism of ${\cal A}$ and $U^*U\in{\cal A}'$.
The aim of the  section is to give a description of the maximal ideal space
$M(\overline{E_*({\cal A}))}$ of the coefficient $C^*$-algebra
$\overline{E_*({\cal A})}=\overline{\{
\bigcup_{n=0}^\infty\delta_*^n({\cal A})\}}$  in terms of the
the maximal ideal space  $M = M({\cal A})$ of
${\cal A}$ and the action of $\delta$ given by (\ref{b-4}).\\

To start with we introduce a number of objects and notation.

Let  $\widetilde{x}\in M(\overline{E_*({\cal A})})$ be a
linear and multiplicative functional on
$\overline{E_*({\cal A})}$. Let us consider  a sequence of
functionals $\xi^n_{\widetilde{x}}:{\cal A}\rightarrow {\bf C}\,$, $n=0,1,...$, defined by the
conditions
\begin{equation}\label{xi2}
\xi^n_{\widetilde{x}}(a)=\delta_*^n(a)(\widetilde{x}),\quad a\in
{\cal A}.
\end{equation}
Since  $\overline{E_*({\cal A})}=\overline{\{
\bigcup_{n=0}^\infty\delta_*^n ({\cal A})\}}$ it follows  that
 the sequence
$\xi^n_{\widetilde{x}}$ determines
$\widetilde{x}$ in a unique way. On the other hand, since  $\delta_*$ is
an endomorphism of $\overline{E_*({\cal A})}$ we have that the
functionals  $\xi^n_{\widetilde{x}}$ are linear and multiplicative
on ${\cal A}$ (may be zero). So either
\begin{equation}
\label{xi0'}
\xi^n_{\widetilde{x}} = x_n \in M ,
\end{equation}
or
\begin{equation}
\label{xi00'}
\xi^n_{\widetilde{x}} = 0 .
\end{equation}
Clearly the mapping
\begin{equation}
\label{xi000'}
{\widetilde{x}} \to (\xi^0_{\widetilde{x}}, \xi^1_{\widetilde{x}}, ...)
\end{equation}
is an injection.

The next theorem is our first step in  the description of
the maximal ideal space
$M(\overline{E_*({\cal A})})$.

\begin{thm}
\label{idealy0}
Let ${\cal A}\subset L(H)$ be a commutative $C^*$-subalgebra,
$1\in{\cal A}$. Let $\delta(a)=UaU^*$ be an endomorphism of
${\cal A}$,   $U^*U\in {\cal A}'$ and
$\alpha:\Delta\rightarrow M({\cal A})$ be a partial mapping
given by  formula (\ref{b-4}) and  $\Delta_n$ be the sets
defined by  (\ref{D_n}), (\ref{D-n}). Then the maximal ideal space
$M(\overline{E_*({\cal A}))}$ of $\overline{E_*({\cal A})}$
is homeomorphic  to a subset of the countable sum of disjoint
 sets. Precisely the mapping (\ref{xi000'}) (along with observations (\ref{xi0'}) and (\ref{xi00'})
defines the topological embedding
\begin{equation}
\label{zawierzgor}
M(\overline{E_*({\cal A})})\hookrightarrow\bigcup_{N=0}^{\infty}M_N\cup
M_\infty,
\end{equation}
where $M_N$ are    the sets of the form
$$
M_N=\{\widetilde{x}=(x_0,x_1,...,x_N,0,...): x_n\in \Delta_n,\,
\alpha(x_{n})=x_{n-1}, \,n=1,...,N\,\},$$
 and   $M_\infty$ is given by
$$
M_\infty=\{\widetilde{x}=(x_0,x_1,...): x_n\in \Delta_n,\,
\alpha(x_{n})=x_{n-1},\, \,n\in {\bf N}\}.
$$
The topology on the sets $M_N$, $N\in{\bf N}\cup\{0\}$, and
$M_\infty$ is induced  by the  base of neighborhoods of  points
$\widetilde{x}\in M_N$
\begin{equation}
\label{otocz1}
O(a_1,...,a_k,\varepsilon)=\{\widetilde{y}\in M_N:
|a_i(x_N)-a_i(y_N)|<\varepsilon, \, \, i=1,...,k\}
 \end{equation}
and respectively  $\widetilde{x}\in M_\infty$
\begin{equation}
\label{otocz2}
 O(a_1,...,a_k,n, \varepsilon)=\{\widetilde{y}\in
\bigcup_{N=n}^{\infty}M_N\cup M_\infty:
|a_i(x_n)-a_i(y_n)|<\varepsilon, \, \, i=1,...,k\}
\end{equation}
 where $\varepsilon > 0$,
$a_i\in A$, and $k,n \in {\bf N}\cup\{0\}$.
\end{thm}
{\bf Proof.} As we have already observed  the mapping (\ref{xi000'}) is an injection.
Considering the righthand part of (\ref{xi000'})
one can come across the following two possibilities.
 \begin{Item}{1)} Let us assume first that some of  the functionals
$\xi^n_{\widetilde{x}}$ are equal to zero. Let $N$ be the least
number such that
$$\xi^{N+1}_{\widetilde{x}}\equiv 0.$$
Note  that for an arbitrary $n\in{\bf N}$ we have
\begin{equation}\label{xi=0}
\xi^{n}_{\widetilde{x}}\neq 0 \Longleftrightarrow
\widetilde{x}(U^{*n}U^{n})=1
\end{equation}
and the family $\{U^{*n}U^n\}_{n\in{\bf N}}$ forms a commutative
decreasing family of projections (see \cite{lebiedodzij}, Proposition 3.5), that is
for $i\leq j$
$$U^{*i}U^iU^{*j}U^j=
U^{*j}U^jU^{*i}U^i= U^{*j}U^j.
$$
Hence, for each $n> N$ we obtain
$$\widetilde{x}(U^{*n}U^n)=\widetilde{x}(U^{*N+1}U^{N+1}U^{*n}U^n)=
\widetilde{x}(U^{*N+1}U^{N+1})\widetilde{x}(U^{*n}U^n)=0,
$$
that is $\xi^{n}_{\widetilde{x}}\equiv 0$ for $n>N$.\\
 Since  for every    $0 \leq
n\leq N$ we have   $\xi^{n}_{\widetilde{x}}\neq 0 $ it follows that  there exists  $x_n\in M(A)$ such that
\begin{equation}\label{x_n}
\xi^n_{\widetilde{x}}(a)=a(x_n), \quad a\in {\cal A},
\end{equation}
and the mapping (\ref{xi000'}) (along with  (\ref{xi0'}) and  (\ref{xi00'})) takes the form
$\widetilde{x} \mapsto (x_0,x_1,...,x_N,0,...)$.\\
Furthermore, since the projections $U^{*i}U^i$, $U^jU^{*j}$
commute we have that
$$U^{*n}U^n U^{*n-1}=U^{*n-1}(U^*U)(U^{n-1} U^{*n-1})=(U^{*n-1}U^{n-1}
U^{*n-1})U^*U=U^{*n}U,$$
$$U^{n-1}U^{*n} U^{n}=(U^{n-1}U^{*n-1})(U^*U)U^{n-1}=U^*U
(U^{n-1}U^{*n-1}U^{n-1})=U^{*}U^{n}.
$$
This along with the fact that $\widetilde{x}(U^{*n}U^n)=1$, for
each $n=1,...,N$, implies that for an arbitrary $a\in {\cal A}$
and $0<n\leq N$ we have
$$a(x_{n-1})=\xi^{n-1}_{\widetilde{x}}(a)=\widetilde{x}(\delta_*^{n-1}(a))=
\widetilde{x}(U^{*n}U^n)\widetilde{x}(\delta_*^{n-1}(a))\widetilde{x}(U^{*n}U^n)=
$$
$$=\widetilde{x}(U^{*n}U^n U^{*n-1}aU^{n-1}U^{*n} U^{n})=\widetilde{x}(U^{*n}UaU^*U^n)=
\widetilde{x}(\delta_*^{n}(\delta(a)))=$$
$$=\xi^{n}_{\widetilde{x}}(\delta(a))=\delta(a)(x_n)=a(\alpha(x_n))$$
where in the  two final equalities we used formulae (\ref{x_n}) and
(\ref{b-4}). Since  ${\cal A}$  separates points of
$M({\cal A})$ we obtain
\begin{equation}\label{alpha(x_n)}
\alpha(x_n)=x_{n-1}, \qquad n\in {\bf N}.
\end{equation}
Hence
we have
\begin{equation}
x_n \in \Delta_n, \qquad 0 \leq n\leq N,
\end{equation}
and thus
$$
\widetilde{x} \mapsto (x_0,x_1,...,x_N,0,...)\in M_N.
$$
\end{Item}

\begin{Item}{2)}
 Now let us assume that for an arbitrary $n\in{\bf N}$ we
have
$$
\xi^n_{\widetilde{x}}\neq 0,
$$
 that is $\xi^n_{\widetilde{x}}\in
M({\cal A})$. Thus, taking $x_n\in M({\cal A})$ such that
(\ref{x_n}) is true we see that  the mapping (\ref{xi000'}) takes the form
$\widetilde{x}\mapsto (x_0,x_1,x_2,...)$.\\
By the same argument as presented in  case 1) we obtain that
\begin{equation}
\alpha(x_n)=x_{n-1}, \qquad  n \geq 1,
\end{equation}
Hence in  view of the  definition of $\Delta_n$, $n\in {\bf Z}$
we have
\begin{equation}
 x_n\in \Delta_n, \qquad n\in {\bf N}.
\end{equation}
And therefore
$$
\widetilde{x} \mapsto (x_0,x_1,x_2,...)\in M_\infty .
$$
\end{Item}
We recall that the mapping  (\ref{xi000'}) is an injection and as we have already proved its right hand
side belongs to $ M_N $ or $M_\infty$ (depending on $\widetilde{x}$). This means that  (\ref{xi000'})
defines an embedding (\ref{zawierzgor}).\\
Finally recall that the topology  $M(\overline{E_*({\cal A})})$ is
 $^*-$weak. That is why a
fundamental system of neighborhoods of a point
$\widetilde{x}=(x_0,x_1,...)\in M_\infty$ is the family of sets of
the form
\begin{equation}\label{4baza}
O(b_1,...,b_k,\varepsilon)=\{\widetilde{y}\in
M(\overline{E_*({\cal A})}):
|b_i(\widetilde{x})-b_i(\widetilde{y})|<\varepsilon, \, \,
i=1,...,k\}
\end{equation}
where $b_i\in \overline{E_*({\cal A})}$, $\varepsilon >0$.
Since
$\overline{E_*({\cal A})}=\overline{\{\bigcup_{n=0}^{\infty}\delta_*^n({\cal A})\}}$
it is enough to take  $b_i=\delta_*^n(a_i)$, $a_i \in
{\cal A}$, $i=0,...,k$, in  (\ref{4baza}) and then
we have that
$$O(b_1,...,b_k,\varepsilon)=\{\widetilde{y}\in
M(\overline{E_*({\cal A})}):
|\delta_*^n(a_i)(\widetilde{x})-\delta_*^n(a_i)(\widetilde{y})|<\varepsilon,
\, \, i=1,...,k\}=$$
$$=\{\widetilde{y}=(y_0,y_1,...)\in
M(\overline{E_*({\cal A})}): |a_i(x_n)-a_i(y_n)|<\varepsilon, \,
\, i=1,...,k\}.$$ If we set
$O(a_1,...,a_k,n,\varepsilon)=O(b_1,...,b_k,\varepsilon)\cap(\bigcup_{N\geq
\, n}M_N\cup M_\infty)$ then we obtain a fundamental system of
neighborhoods
from the thesis.\\
In the case of a point  $\widetilde{x}=(x_0,x_1,...,x_N,0,...)\in M_N$
 we set
$$
O(a_1,...,a_k,\varepsilon)=O(a_1,...,a_k,N,\varepsilon)\cap
 M_N.
$$
Thus  formulae (\ref{otocz1}) and (\ref{otocz2}) define a base of
neighborhoods of a point
 $\widetilde{x}\in M(\overline{E_*({\cal A})}) \hookrightarrow\bigcup_{N=0}^{\infty}M_N\cup
M_\infty$ and the proof is complete. \\
\begin{rem}
\label{ostatniauwaga}
\em
A very useful fact is that the set of operators of the form
$b=a_0+\delta_*(a_1)+...+\delta_*^N(a_N)$, where
$a_0,a_1,...,a_N\in A$,  is  dense in
$\overline{E_*({\cal A})}$ (see \cite{lebiedodzij}, Proposition  3.7.)
and by definition  (\ref{x_n}) of
functionals creating the sequence  $\widetilde{x}=(x_0,...)$ we
have that
$$
 b(\widetilde{x})=a_0(x_0)+a_1(x_1)+...+a_N(x_N),
$$
when
$
\widetilde{x}=(x_0,...)\in \bigcup_{n=N}^\infty M_n\cup
M_\infty $\\
 and
$$
     b(\widetilde{x})=a_0(x_0)+a_1(x_1)+...+a_n(x_k),
$$
when $\widetilde{x}=(x_0,...,x_k,0,...)\in \bigcup_{n=0}^N M_n.$
\end{rem}

The theorem just proved can be considered us an 'upper estimate'
for $M(\overline{E_*({\cal A})})$. The next theorem
\ref{idealy0.0.5} presents the corresponding 'lower estimate'.
Before passing to
 its statement let us stress  that the previous theorem tells us that any element
$\widetilde{x}\in M(\overline{E_*({\cal A})}) $  generates (defines uniquely) either an element of $M_N $
or an element of $M_\infty$, on the other hand given a  sequence $(x_0,x_1,...,x_N,0,...) \in M_N $
or $(x_0,x_1,...,x_n,...) \in M_\infty$ one cannot say in advance that it is generated by a certain
element $\widetilde{x}\in M(\overline{E_*({\cal A})}) $ (this will have place only if this sequence is of the form
(\ref{xi000'})). Theorem \ref{idealy0.0.5} tells that all the sequences from the  subset $\widehat{M}_N\subset M_N$
(which in general can be {\em essentially} smaller than $M_N$) and all the sequences from $M_\infty$ are really generated
by certain elements $\widetilde{x}\in M(\overline{E_*({\cal A})}) $.
\begin{thm}
\label{idealy0.0.5} Let ${\cal A}\subset L(H)$ be a commutative
$C^*$-algebra, $1\in{\cal A}$, and let $U\in L(H)$ be  a partial
isometry such that $U^*U\in{\cal A}'$. In view of theorem
\ref{idealy0} we can treat $M(\overline{E_*({\cal A})})$ as a
subset of $\bigcup_{N\geq \, 0}M_N\cup M_\infty$ on the other hand
the following inclusion holds
\begin{equation}
\label{zawierzdol}
\bigcup_{N=0}^{\infty}\widehat{M}_N\cup M_\infty \subset
M(\overline{E_*({\cal A})})
\end{equation}
where $\widehat{M}_N$ are the sets of the form
$$
    \widehat{M}_N=\{\widetilde{x}=(x_0,x_1,...,x_N,0,...): x_N\in
    \Delta_N,\,x_N\notin \Delta_{-1},\, \alpha(x_{n})=x_{n-1}\}
$$
and $ M_\infty=\{\widetilde{x}=(x_0,x_1,...): x_n\in \Delta_n,\,
\alpha(x_{n})=x_{n-1},\, \,n\geq 1\}.$
\end{thm}
\textbf{Proof}. Let us  prove first that $M_\infty \subset
M(\overline{E_*({\cal A})})$. Let $(x_0,x_1,...)$ be an
arbitrary sequence of elements  of  $M({\cal A})$ satisfying
the condition \ \  $\alpha(x_{n})=x_{n-1}$ \ \ for each \ \ $n\geq 1$.
We shall show that there exists a  linear and
multiplicative functional $\widetilde{x}\in
M(\overline{E_*({\cal A})})$ that generates this sequence
(that is all the equations (\ref{x_n}), $n=0,1,2, ... $ where $ \xi^n_{\widetilde{x}}$ are  given by (\ref{xi2}) are satisfied).\\
Indeed, let us consider the following sets
\begin{equation}\label{X^n}
  \widetilde{X}^{n}=\{\widetilde{x}\in M(\overline{E_*({\cal A})}):\,
  \forall_{a\in A}\ \, \xi^n_{\widetilde{x}}(a)=a(x_n)\,\},
  \qquad  n=0,1...\,
\end{equation}
where $ \xi^n_{\widetilde{x}}$ are given by (\ref{xi2}). We have
that

\begin{description}
\item[1)] $\widetilde{X}^{n}\neq \emptyset$. This is due to the fact that  (\ref{X^n})
is the set of all extensions of a certain multiplicative
functional from  $\delta_*^{n}({\cal A})$ up to
$\overline{E_*({\cal A})}$. It is known that for any
multiplicative functional on a commutative $C^*$-subalgebra there
exists an extension up to a multiplicative functional on a
larger commutative  $C^*$-algebra.

\item[2)] $\widetilde{X}^{n}$ is closed. This  follows  from the
definition of weak$^*$ convergence.

\item[3)] $\widetilde{X}^{0}\supset\widetilde{X}^{1}\supset...\supset\widetilde{X}^{n}\supset...\,
\,.$  Indeed, if $\widetilde{x}\in \widetilde{X}^{n}$ then
$$
   a(x_{n-1})=a(\alpha(x_n))=\delta(a)(x_n)= \xi^n_{\widetilde{x}}(\delta(a))=
   \xi^{n-1}_{\widetilde{x}}(a),\qquad a\in {\cal A},
$$
that is $\widetilde{x}\in \widetilde{X}^{n-1}$.
\end{description}
Hence the family of sets $\widetilde{X}^{n}$ forms a decreasing
sequence of nonempty compact sets and thus
$\bigcap_{n=0}^{\infty}\widetilde{X}^{n}\neq \emptyset$. From  the definition  of $\widetilde{X}^{n}$
it follows that any point
$\widetilde{x} \in \bigcap_{n=0}^{\infty}\widetilde{X}^{n}$ generates the same  given sequence
$$
( \xi^0_{\widetilde{x}},  \xi^1_{\widetilde{x}},  \xi^2_{\widetilde{x}}, ... ) = (x_0,x_1,x_2,... ).
$$
  But we have already observed that the mapping
$$
\widetilde{x} \to  ( \xi^0_{\widetilde{x}},  \xi^1_{\widetilde{x}},  \xi^2_{\widetilde{x}}, ... )
$$
is injective and therefore
$$
 \bigcap_{n=0}^{\infty}\widetilde{X}^{n} =\{\widetilde{x}\}
$$
consists  of a single point.\ \\

Now, let  $(x_0,x_1,...,x_N,0,...)\in\widehat{M}_N$. Let us
consider the sets $\widetilde{X}^{n}$ defined by  (\ref{X^n}), but
only for $n=0,1,...,N$. The above argument shows that these sets
are decreasing and nonempty. Let $\widetilde{x}
\in\widetilde{X}^{N}$. To identify  $\widetilde{x}$ with the
sequence $( \xi^0_{\widetilde{x}},  \xi^1_{\widetilde{x}},  ... , \xi^N_{\widetilde{x}},0, ... )=(x_0,x_1,...,x_N,0,...)$ given by (\ref{xi000'}) it
is enough to show that
$$
        \xi^{N+1}_{\widetilde{x}}\equiv 0.
$$
Let us assume the opposite, that is $\xi^{N+1}_{\widetilde{x}}\neq
0$. Then from the first part of the proof of theorem \ref{idealy0}
 we obtain that
$\widetilde{x}=(x_0,x_1,...,x_N,x_{N+1},...)$ where
$\alpha(x_{N+1})=x_{N}$, which  contradicts  the fact that
$x_N\notin \Delta_{-1}=\alpha(\Delta_1)$. The proof is complete. \\

 In general in the preceding theorems \ref{idealy0}, \ref{idealy0.0.5}
  none of the
inclusions  (\ref{zawierzgor}) or (\ref{zawierzdol}) can be replaced
by the equality. However,  as the next result tells in the situation when ${\cal A}$ contains
 $\delta_*(1)=U^*U$ we have
  the equality in  formula
(\ref{zawierzdol}) and therefore  it gives  the full description of the maximal ideal space of
$M(\overline{E_*({\cal A})})$.

\begin{thm}\label{idealy0.1}
Let $A\subset L(H)$ be commutative  $C^*$-algebra containing the
identity. Let $\delta(\cdot)=U(\cdot )U^*$ be an endomorphism of
${\cal A}$. Moreover, let
$$
U^*U\in {\cal A}.
$$
Then the maximal ideal space  $M(\overline{E_*({\cal A})})$
 of the algebra $\overline{E_*({\cal A})}$ is homeomorphic to
 the countable sum of disjoint closed and open sets $\widehat{M}_N$ and
the
 closed set $M_\infty$ (we admit empty sets)
\begin{equation}
\label{equality}
M(\overline{E_*({\cal A})})=\bigcup_{N=0}^{\infty}\widehat{M}_N\cup
M_\infty
\end{equation}
where
$$
\widehat{M}_N=\{\widetilde{x}=(x_0,x_1,...,x_N,0,...): x_N\in
\Delta_N,\,x_N\notin \Delta_{-1},\, \alpha(x_{n})=x_{n-1}\},$$
$$
M_\infty=\{\widetilde{x}=(x_0,x_1,...): x_n\in
\Delta_n\cap\Delta_{-\infty},\, \alpha(x_{n})=x_{n-1},\,
\,n\geq 1\}.
$$
The topology on $\bigcup_{N=0}^{\infty}\widehat{M}_N\cup M_\infty$
is defined  by a fundamental system of neighborhoods of  points
$\widetilde{x}\in \widehat{M}_N$
\begin{equation}\label{otocz1.1}
O(a_1,...,a_k,\varepsilon)=\{\widetilde{y}\in \widehat{M}_N:
|a_i(x_N)-a_i(y_N)|<\varepsilon, \, \, i=1,...,k\}
 \end{equation}
and respectively   $\widetilde{x}\in M_\infty$
$$ O(a_1,...,a_k,n, \varepsilon)=\{\widetilde{y}\in
\bigcup_{N=n}^{\infty}\widehat{M}_N\cup M_\infty:
|a_i(x_n)-a_i(y_n)|<\varepsilon, \, \, i=1,...,k\}
$$
 where $\varepsilon > 0$,
$a_i\in A$ and $k,n \in {\bf N}\cup\{0\}$.
\end{thm}
{\bf Proof.} In view of  theorems  \ref{idealy0} and
\ref{idealy0.0.5} to prove the equality (\ref{equality}) it is enough to prove that for an arbitrary
$\widetilde{x}\in M(\overline{E_*(A)})$ we have
$$
\widetilde{x}\in M_N\Longrightarrow\widetilde{x}\in
\widehat{M}_N.
$$
 Suppose on the contrary that $\widetilde{x}\in M_N$ and  $\widetilde{x}\notin
\widehat{M}_N$. Then $\widetilde{x}=(x_0,x_1,...,x_N,0,...)$,
$x_n\in M({\cal A})$ and there exists $x_{N+1}\in
\Delta_1\subset M({\cal A})$ such that $\alpha(x_{N+1})=x_N$.
By definition (\ref{x_n}) of functionals $x_n$ and definition
 (\ref{xi2}) of functionals $\xi_{\widetilde{x}} ^n$
we have that  $\widetilde{x}(U^{*n}aU^n)=a(x_n)$, for each $a\in
A$ and $n=0,...,N$, and $\widetilde{x}(U^{*n}aU^n)=0$ when $n> N$.
In particular we have that
$$
\widetilde{x}(U^{*N}U^N)=1 \quad {\rm and} \quad
\widetilde{x}(U^{*N+1}U^{N+1})=0.
$$
By  formula  (\ref{b-4})  we obtain that
$$
\widetilde{x}(U^{*N}aU^N)=a(x_N)=a(\alpha(x_{N+1}))=\delta(a)(x_{N+1}),\qquad
a\in A.
$$
Setting  $a=U^*U\in A$ in the above formula we have that
$\delta(U^*U)(x_{N+1})=0$. On the other hand for
$\delta(U^*U)(x_{N+1})$ we have
 $$
\delta(U^*U)(x_{N+1})=x_{N+1}(UU^*UU^*)=x_{N+1}(UU^*)\cdot
x_{N+1}(UU^*)=1.
 $$
Thus we arrived at a  contradiction  and the proof of (\ref{equality})  is complete.\\

The closedness and openness of the sets $\widehat{M}_N$ follows from the relation
$$
\widetilde{x}\in M_N\Longleftrightarrow\widetilde{x}\in
\widehat{M}_N \Longleftrightarrow \{ \widetilde{x}(U^{*N}U^N)=1 \quad {\rm and} \quad
\widetilde{x}(U^{*N+1}U^{N+1})=0 \}.
$$
Further the openness of the sets $\widehat{M}_N$ implies the closedness of $M_\infty$
which finishes the proof of the theorem.

\begin{rem}\label{algA1}
\emph{ In fact the theorem just proved gives us the key to obtain the complete
description of  $
M(\overline{E_*({\cal A})})
$
in the general situation. Indeed, if  $U^*U\notin {\cal A}$ then  one can consider the
$C^*$-algebra ${\cal A}_1=\langle {\cal A},U^*U \rangle$
generated by ${\cal A}$ and $U^*U$. Since
$$
\delta(U^*U)=UU^*UU^*=\delta(1)\delta(1)\in {\cal A}
$$
we have that $\delta:{\cal A}_1\rightarrow {\cal A}$. By applying
the foregoing theorem to the algebra ${\cal A}_1$ and the operator
$U$ one obtains the full description of $M(\overline{E_*({\cal
A}_1)})$. However, as
$$
\overline{E_*({\cal A}_1)}=\overline{E_*({\cal A})}
$$
 we
have
$$
M(\overline{E_*({\cal A}_1)})=M(\overline{E_*({\cal A})}).
$$}
\end{rem}

The preceding results often can be improved (simplified) when we know some
specific  features of a partial isometry $U$. For instance, if $U$
is an isometry, i.e. $U^*U=1$, then by theorem \ref{takietamtwierdzenie}
the mapping $\alpha$ generated
by the endomorphism $\delta$ is surjective   and hence all the sets
$\widehat{M}_N$ are empty. Thus by theorem \ref{idealy0.1} we
obtain the following

\begin{corl}\label{idealy1}
Let $A\subset L(H)$ be a commutative $C^*$-algebra containing the
identity and  $U$ be an isometry such that  $\delta(\cdot)=U(\cdot)U^*$
is an endomorphism of ${\cal A}$. Then  the mapping $\alpha:\Delta\rightarrow
M({\cal A})$ defined by (\ref{b-4}) is a  surjection  and the
spectrum  $M(\overline{E_*({\cal A})})$ of $\overline{E_*(A)}$
has  the form
\begin{equation}
M(\overline{E_*({\cal A})})=\{\widetilde{x}=(x_0,x_1,...):
x_n\in \Delta_n,\, \alpha(x_{n+1})=x_{n},\, \,n\geq 0\}
\end{equation}
A fundamental system of neighborhoods of a point $\widetilde{x}$
is family of sets
$$
O(a_1,...,a_k,n,\varepsilon)=\{\widetilde{y}\in
M(\overline{E_*(A)}): |a_i(x_n)-a_i(y_n)|<\varepsilon, \, \,
i=1,...,k\}
$$
where $n\in {\bf N}\cup\{0\}$, $\varepsilon > 0$ and  $a_i\in
{\cal A}$, $ 1 \leq i\leq k$.
\end{corl}

\end{document}